\documentclass[hidelinks,onefignum,onetabnum]{siamart220329}

\usepackage{lipsum}
\usepackage{amsfonts}
\usepackage{graphicx}
\usepackage{epstopdf}
\usepackage{algorithmic}
\ifpdf
  \DeclareGraphicsExtensions{.eps,.pdf,.png,.jpg}
\else
  \DeclareGraphicsExtensions{.eps}
\fi
\usepackage{amsopn}
\usepackage{makecell}
\usepackage{siunitx}
\usepackage{colortbl}
\usepackage{multirow}
\newsiamremark{remark}{Remark}
\newsiamremark{hypothesis}{Hypothesis}
\crefname{hypothesis}{Hypothesis}{Hypotheses}
\newsiamthm{claim}{Claim}

\headers{Physics-informed NVAR}{J. H. Adler, S. Hocking, X. Hu, and S. Islam}

\title{Physics-informed nonlinear vector autoregressive models for the prediction of dynamical systems\thanks{Submitted to the editors \today.
\funding{This work was funded by Tufts University DISC Seed Funding: On the Predictability Limits of Nonlinear Chaotic Systems: Can Deep Learning Methods Help?}}}

\author{
    James H. Adler\thanks{Department of Mathematics, Tufts University, Medford, MA (\email{james.adler@tufts.edu}, \email{samuel.hocking@tufts.edu}, \email{xiaozhe.hu@tufts.edu}).}
    \and Samuel Hocking\footnotemark[2]
    \and Xiaozhe Hu\footnotemark[2]
    \and Shafiqul Islam\thanks{Water Diplomacy, Civil and Environmental Engineering, and The Fletcher School of Law and Diplomacy, Tufts University, Medford, MA (\email{shafiqul.islam@tufts.edu}).}
}

\usepackage{amsopn}

\ifpdf
\hypersetup{
  pdftitle={Physics-informed nonlinear vector autoregressive models for the prediction of dynamical systems},
  pdfauthor={J. H. Adler, S. Hocking, X. Hu, and S. Islam}
}
\fi

\begin{document}

\maketitle

% REQUIRED
\begin{abstract}
Machine learning techniques have recently been of great interest for solving differential equations. Training these models is classically a data-fitting task, but knowledge of the expression of the differential equation can be used to supplement the training objective, leading to the development of physics-informed scientific machine learning. In this article, we focus on one class of models called nonlinear vector autoregression (NVAR) to solve ordinary differential equations (ODEs). Motivated by connections to numerical integration and physics-informed neural networks, we explicitly derive the physics-informed NVAR (piNVAR) which enforces the right-hand side of the underlying differential equation regardless of NVAR construction. Because NVAR and piNVAR completely share their learned parameters, we propose an augmented procedure to jointly train the two models. Then, using both data-driven and ODE-driven metrics, we evaluate the ability of the piNVAR model to predict solutions to various ODE systems, such as the undamped spring, 
a Lotka-Volterra predator-prey nonlinear model, and the chaotic Lorenz system.
\end{abstract}

% REQUIRED

% REQUIRED

\section{Introduction}\label{sec:intro}

Ordinary differential equations (ODEs) are used across domains in science and engineering to model phenomena and systems of interest including simple models of motion, springs, population growth, resource utilization, and more (see, e.g. \cite{blanchard_differential_2006, borgers_introduction_2017}). We consider a prototypical initial value problem of the following general form,
\begin{align}\label{eqn:IVP}
    \frac{d\mathbf{x}}{dt}=\mathbf{x}'=\mathbf{f}(\mathbf{x}),\quad \mathbf{x}(0)=\mathbf{x}_0,
\end{align}
with $\mathbf{x}\in\mathbb{R}^d$, $\mathbf{f}:\mathbb{R}^d \mapsto \mathbb{R}^d$. For certain $\mathbf{f}$, a closed-form solution may be obtained. Otherwise, we must turn to numerical techniques to obtain a solution. One ubiquitious class of methods are numerical integration schemes such as linear multistep methods (LMMs) and Runge-Kutta methods (RKMs). Such methods are cheap to evaluate and are supported by an ample body of stability, consistency, and convergence theory. However, they require confident knowledge of $\mathbf{f}$ and have no data-driven mechanism to operate in the absence of such knowledge. On the other hand, recurrent neural network (RNN) architectures, which originated as tools for data-driven timeseries forecasting \cite{bianchi_recurrent_2017}, offer a counterpoint to ODE-centric numerical methods. Long short-term memory (LSTM) \cite{hochreiter_long_1997} and gated recurrent unit (GRU) \cite{cho_learning_2014} models were meant to alleviate problems with vanishing or exploding gradients seen in densely connected feed-forward neural networks and designed to produce a sequence of points beginning from a short sequence of previous observations by ``learning" latent patterns in the training data. In contrast to LMMs and RKMs, RNNs are comparatively costly to train and evaluate.

A third class of numerical methods for ODEs is comprised of so-called reservoir computers, named for their constructed ``reservoir" of dynamics available to model the actual dynamics generated and described by the training data. Echo state networks (ESNs) \cite{jaeger_echo_2001, jaeger_optimization_2007, yildiz_re-visiting_2012} feature randomly initialized input weights, a random and sparsely connected internal network, and trained output weights. Though ESNs are fundamentally data-driven models, variants of ESN which utilize knowledge of (\ref{eqn:IVP}) have also been developed \cite{pathak_hybrid_2018, rodrigues_physics-informed_2019, doan_learning_2020}. ESN training is accomplished by linear least-squares and is therefore cheap and easy to implement. However, the non-deterministic model initialization poses challenges to analysis and robust understanding of the large parameter space. Special cases of ESNs have been analyzed and have led to the development of deterministic nonlinear vector autoregression (NVAR) models \cite{bollt_explaining_2021, gauthier_next_2021}. NVAR replaces the randomly initialized internal network with an intentionally chosen non-random nonlinear function while retaining the same least-squares training process. As a result, we chose to focus on the NVAR model type for this work.

RNNs and reservoir computers are conventionally data-driven models. That is, their respective training processes, with the exception of any regularization or non-overfitting aspects, focus entirely on fitting the coordinate $\mathbf{x}$ with the training data. No direct consideration is given to the satisfaction of the physical constraint, $\mathbf{f}(\mathbf{x})$. Physics-informed neural networks (PINNs) were recently developed \cite{raissi_physics-informed_2019} with a focus on solving partial differential equations (PDEs) (see, e.g. \cite{cai_physics-informed_2021, cai_physics-informed_2021-1, cuomo_scientific_2022, mao_physics-informed_2020}), in which the primary neural network enforces initial and boundary conditions, and a secondary physics-informed network enforces the statement of the PDE. DeepXDE \cite{lu_deepxde_2021}, its implementation, uses gradient-based training, which is relatively slow and computationally costly. By comparison, a physics-informed reservoir computer should be easy to explicitly derive and cheap to train in a linear least-squares process. 

A physics-informed ESN was first proposed and demonstrated in \cite{doan_learning_2020}. However, the model training process is no longer completely linear since it includes a second general optimization state solved by L-BFGS-B. Automatic-differentiated physics-informed ESNs (API-ESNs) \cite{racca_automatic-differentiated_2021} built upon this work by restoring the linear training process through the explicit derivation of the model's time derivative. While the API-ESN successfully implements the physics-informed concepts from PINNs and is cheap to train, the model retains the attributes that challenge conventional ESNs: random model initialization, numerous hyperparameters, and lack of transparency and explainability. We further advance the concept of a physics-informed general-purpose reservoir computer by proposing physics-informed NVAR (piNVAR). Similar to API-ESN, piNVAR is obtained by the explicit derivation of the time derivative of the NVAR output through the chain rule and training remains completely linear. However, piNVAR retains the same advantages as NVAR: deterministic construction, few hyperparameters, and flexible yet explainable nonlinearity. These positive attributes facilitate our development of strong theoretical understanding and broad numerical examination. Using the simple expressions of NVAR and piNVAR update formulae, we show the coupling between NVAR and piNVAR and firmly connect our findings to the PINN theory developed in \cite{raissi_physics-informed_2019}. We demonstrate that piNVAR, through completely shared parameters with NVAR, propogates an update in the time derivative of the model prediction which conforms to the governing ODE. Most importantly, using piNVAR requires only easily expressed gradient information and incremental computation, so is advantageous for all NVAR model constructions. Through a broad cross-validation testing routine across multiple test problems, a grid of model parameters, and multiple explainable nonlinear state functions, we demonstrate that physics-informed NVAR training yields substantial improvements to both data-driven and physics-informed evaluation metrics.

In Section \ref{sec-nvar}, we provide an overview of the NVAR model and review the linear least-squares training process. Section \ref{sec-pinvar} begins with a brief synopsis of PINNs, followed by the derivation of the piNVAR update formula and a discussion of the relationship between the NVAR and piNVAR models. Augmented training procedures are also given. Then, we demonstrate the effectiveness of piNVAR prediction measured by both data-driven and ODE-driven metrics in Section \ref{sec-num-exp}. We conclude with a synthesis of our contributions, findings, and thoughts on future directions in Section \ref{sec-conclusions}.

\section{Nonlinear vector autoregression (NVAR) models}\label{sec-nvar}

Before describing the construction and training of an NVAR model, we fix some notation. Let $\mathbf{x}(t)\in\mathbb{R}^{d}$ be a point at time $t$ which solves the $d$-dimensional initial value problem (\ref{eqn:IVP}). In this case, $\mathbf{x}(t)$ represents the system's coordinates in state space. We also refer to such coordinates as a point. Let time $t\geq 0$ be discretized by a fixed time step $h=\Delta t>0$ and denote $t_{k}=t_0+kh$ where it is generally assumed that $t_0=0$. The notations $\mathbf{x}(t_k)$ and $\mathbf{x}(t_0+kh)$ are equivalent, while $\mathbf{x}_k$ represents a point which approximately solves the IVP at time $t=t_k$. Finally, $x_i(t_k)$ and $x_{k,i}$ denote the $i^{th}$ entry of the vectors $\mathbf{x}(t_k)$ and $\mathbf{x}_{k}$, respectively.

The central mechanism of NVAR \cite{gauthier_next_2021,shahi_prediction_2022} is the prediction of $\mathbf{x}_{k+1}$ by a linear combination of nonlinear functions of the current iterate, $\mathbf{x}_k$, and $p-1$ previous iterates. The model inputs are controlled by two parameters. Lookback, $p \geq 1$, specifies the total number of utilized data points. Sampling frequency, $s \geq 1$, sets the interval at which the utilized data points are drawn. Given $p$ and $s$, we denote the corresponding input vector at time $t=t_k$ by
\begin{align*}
    \mathbf{y}_k^{p,s}= \begin{bmatrix}
        \mathbf{x}_{k} \\ \mathbf{x}_{k-s} \\ \vdots \\ \mathbf{x}_{k-(p-1)s}
    \end{bmatrix} \in \mathbb{R}^{pd}.
\end{align*}
In addition to the input parameters, the primary degree of freedom in NVAR construction is the choice of nonlinear state function, $\mathbf{h}:\mathbb{R}^{pd}\rightarrow\mathbb{R}^m$. Here, $m$ depends on the particular choice of $\mathbf{h}$ and we typically have $m > d$. The state function serves an analagous purpose to the activation function in neural networks. Here, however, the trained weights matrix, $\mathbf{W}\in\mathbb{R}^{d \times m}$, acts on the evaluated state function, rather than the activation function acting on the affine linear transformation from a particular neuron. NVAR's linear dependence on the trained parameters arises from this subtle difference. In this work, $\mathbf{W}$ is trained to facilitate the update,
\begin{align}\label{int-NVAR-pred}
    \mathbf{x}_{k+1} &= \mathbf{x}_{k} + \mathbf{W}\mathbf{h}(\mathbf{y}_k^{p,s}),
\end{align}
and we denote the target of the NVAR linear combination by,
\begin{align}\label{z-nvar-tgt}
    \mathbf{z}_k &= \mathbf{W}\mathbf{h}(\mathbf{y}_k^{p,s}) = \mathbf{x}_{k+1} - \mathbf{x}_{k}.
\end{align}
Here, we regard $\mathbf{W}\mathbf{h}(\mathbf{y}_k^{p,s})$ as an approximate definite integral such that
\begin{align}\label{I-approx}
    \mathbf{W}\mathbf{h}(\mathbf{y}_k^{p,s}) \approx \int_{kh}^{(k+1)h} \mathbf{f}(\mathbf{x}(t))dt, 
\end{align}
and hence call Equation \cref{int-NVAR-pred} the \textit{integration update formula}.

Suppose $\mathcal{D}=\{\mathbf{x}(t_1),\ldots,\mathbf{x}(t_N) \}$ is the available dataset of $N$ data points, indexed from 1 to $N$, and equally spaced by time step $\Delta t$. To train a model, we must determine a suitable subset of training points. Given model parameters $p$ and $s$, allow at least $a>(p-1)s$ excess points before marking the beginning of the training set. Choose a suitable number of training points, $T$, and let the training indices be $a,a+1,\ldots,a+T-1$. Populate the training target matrix $\mathbf{Z}_{tr}\in\mathbb{R}^{d \times T}$ with each column equal to the single time-step target $\mathbf{z}(t_{a+k})=\mathbf{x}(t_{a+k+1})-\mathbf{x}(t_{a+k})$ for the integration update formula (\ref{int-NVAR-pred}), for $k=0,\ldots,T-1$.
\begin{align}
    \mathbf{Z}_{tr} = \begin{bmatrix}
			\mathbf{z}(t_{a}) & \ldots & \mathbf{z}(t_{a+T-1})
		\end{bmatrix}
\end{align}
Then, the training state matrix $\mathbf{H}_{tr}\in \mathbb{R}^{m\times T}$ is formed with columns
$\mathbf{h}(\mathbf{y}_{a+k}^{p,s})$, $k=0,\ldots,T-1$.
\begin{align}
    \mathbf{H}_{tr}=\begin{bmatrix}
			\mathbf{h}(\mathbf{y}_a^{p,s}) & \ldots & \mathbf{h}(\mathbf{y}_{a+T-1}^{p,s})
		\end{bmatrix}
\end{align}
Training $\mathbf{W}$ is accomplished by minimizing the following objective function over all matrices $\mathbf{W}\in\mathbb{R}^{d\times m}$,
\begin{align}\label{orig_obj}
	\min_{\mathbf{W}\in\mathbb{R}^{d\times m}} \quad w_d g_d(\mathbf{W})  + r g_r(\mathbf{W}).
\end{align}
In the above, we define the data-loss term by $g_d(\mathbf{W}):=\lVert \mathbf{W}\mathbf{H}_{tr} - \mathbf{Z}_{tr} \rVert_F^2$ and the regularization term by $ g_r(\mathbf{W}):= \lVert \mathbf{W} \rVert_F^2$.
Here, $w_d,r\geq 0$ are the data weight and regularization parameters, respectively, and $\lVert\ \cdot \ \rVert_F$ denotes the Frobenius norm. This objective is equivalent to
\begin{align}\label{equiv_obj}
	\min_{\mathbf{W}\in\mathbb{R}^{d\times m}} \quad \big\lVert \mathbf{L}(\mathbf{W})\big\rVert_F^2,
\end{align}
where
\begin{align}\label{equiv_L}
    \mathbf{L}(\mathbf{W}) = \mathbf{W}\begin{bmatrix}
    \sqrt{w_d}\ \mathbf{H} \\
    \sqrt{r}\ \mathbf{I}
    \end{bmatrix} - \begin{bmatrix}
    \sqrt{w_d}\ \mathbf{Z}_{tr} \\
    \mathbf{0}
    \end{bmatrix}.
\end{align}
The optimization problem (\ref{equiv_obj}) is solved by linear least squares and is thus quite straightforward either explicitly or by a range of off-the-shelf implementations.

\section{Physics-informed NVAR}\label{sec-pinvar}

In Section \ref{sec-nvar}, we reviewed the generic NVAR model and highlighted the substantial flexibility embedded in the user-defined state function, $\mathbf{h}$. Some informed choices of the state function are discussed in further detail in \cite{HockingSamuelC2024PoDS}. In some cases, such as the linear multistep method and Runge-Kutta representations by NVAR, the state function includes evaluations of the right-hand side of the underlying differential equation, $\mathbf{f}$. However, given an arbitrary state function which does not include evaluations of $\mathbf{f}$ and which may not effectively preserve the underlying dynamics, we seek a method through which we can enforce those dynamics. In this Section, we derive the physics-informed NVAR which completely shares its trained parameters with the data-driven NVAR, and demonstrate how both models are simultaneously trained through the addition of an ODE-fit term to the conventional NVAR training objective. 

\subsection{PINNs}\label{subsec-pinns}

While the conventional NVAR model is data-driven by virtue of purely data-driven training, we take inspiration from the incorporation of knowledge of the underlying differential equation to data-driven black-box feed-forward neural networks which resulted in PINNs. Consider a feed-forward neural network as in \cite{lu_deepxde_2021, zeinhofer_unified_2024} and the general form of a PDE parameterized by $\boldsymbol{\lambda}$,
\begin{align}\label{PDE-form}
    f\bigg(\mathbf{x},\frac{\partial u}{\partial x_1},\ldots,\frac{\partial u}{\partial x_d},\frac{\partial^2 u}{\partial x_1\partial x_1},\ldots,\frac{\partial^2 u}{\partial x_1\partial x_d},\ldots,\boldsymbol{\lambda} \bigg)&=0,\quad \mathbf{x}\in\Omega, \\
    \mathcal{B}(u,\mathbf{x})&=0,\quad \mathbf{x}\in\partial\Omega,
\end{align}
where $\mathbf{x}=[x_1 \ \ldots \ x_d]^T$, $\Omega \subset \mathbb{R}^d$, time is included, if applicable, as a dimension of $\mathbf{x}$, and $\mathcal{B}(u,\mathbf{x})$ are general initial and boundary conditions \cite{lu_deepxde_2021}. The dimension of $u$ is left arbitary,
\begin{align}
    u: \mathbf{R}^d\rightarrow \mathbf{R}^m,\quad m \geq 1.
\end{align}
A PINN solution to (\ref{PDE-form}) should approximate $u$ by a feed-forward neural network $u^{\theta}: \mathbf{R}^d\rightarrow \mathbf{R}^m$, parameterized by $\theta$. Borrowing some notation from \cite{zeinhofer_unified_2024}, we let $L \in \mathbb{N}$ denote the number of layers or depth, $N_l\in\mathbb{N}$, $l=0,\ldots,L$ denote number of units or width of layer $l$, with the restriction that $N_0=d$ and $N_L=m$. Finally, let each $(\mathbf{W}_l,\mathbf{b}_l)$, $l=1,\ldots,L$, with conforming dimensions $\mathbf{W}\in\mathbb{R}^{N_l\times N_{l-1}}$ and $\mathbf{b}\in\mathbb{R}^{N_l}$, be a matrix-vector pair that determines the affine linear map,
\begin{align}
    \mathbf{T}_l&:\mathbb{R}^{N_{l-1}}\rightarrow\mathbb{R}^{N_{l}}, \\
    \mathbf{T}_l(\mathbf{y}) &\mapsto \mathbf{W}_l\mathbf{y}+\mathbf{b}_l.
\end{align}
Then, $\theta = \big((\mathbf{W}_1,\mathbf{b}_1),\ldots, (\mathbf{W}_L,\mathbf{b}_L)\big)$, and we call the parameter space $\Theta$ the space of all conforming parameter configurations $\theta$. Let $\sigma:\mathbb{R}\rightarrow\mathbb{R}$ be a nonlinear activation function that is applied entry-wise to $\mathbf{T}_l(\mathbf{y})$, and $u_\theta$ is defined recursively as
\begin{align}
    u^\theta_0(\mathbf{x}) &= \mathbf{x} \in \mathbb{R}^d \\
    u^\theta_l(\mathbf{x}) &= \sigma(\mathbf{W}_l u^\theta_{l-1}(\mathbf{x})+\mathbf{b}_l) \in \mathbb{R}^{N_l},\quad  l=1,\ldots,L-1 \\
    u^\theta_L(\mathbf{x}) &= \mathbf{W}_L u^\theta_{L-1}(\mathbf{x})+\mathbf{b}_L \in \mathbb{R}^{N_m}.
\end{align}
The action of $u^\theta$ follows from the successive application of each layer,
\begin{align}
    u^\theta(\mathbf{x}) &= u^\theta_L(u^\theta_{L-1}(\ldots u^\theta_1(\mathbf{x}))).
\end{align}
In \cite{raissi_physics-informed_2019}, the authors take the helpful abstract view that the neural network $u^\theta$ should approximate the initial and boundary data while a secondary neural network with shared parameters, $f^\theta$, derived from $u^\theta$ by automatic differentiation, perhaps using techniques described in \cite{baydin_automatic_2018}, should approximate and enforce $f$ on $\Omega$. Since the parameters, $\theta$, are shared between $u^\theta$ and $f^\theta$, they can be learned by minimizing a joint objective function,
\begin{align}
    \mathcal{L}(\theta, \mathcal{T}_b,\mathcal{T}_f) &= w_b\underbrace{\mathcal{L}_b(\theta, \mathcal{T}_b)}_{\text{boundary loss}}+\ w_f\underbrace{\mathcal{L}_f(\theta, \mathcal{T}_f)}_{\text{PDE loss}},
\end{align}
where
\begin{align}
    \mathcal{T}_b&=\{\mathbf{x}_1,\ldots,\mathbf{x}_{T_b} \}\subseteq \partial\Omega, \\
    \mathcal{T}_f&=\{\mathbf{x}_1,\ldots,\mathbf{x}_{T_f} \}\subseteq \Omega, \\
    \mathcal{L}_b(\theta, \mathcal{T}_b) &= \frac{1}{T_b}\sum_{\mathbf{x}\in\mathcal{T}_b} \lVert \mathcal{B}(u^\theta, \mathbf{x})\rVert_2^2,\\
    \mathcal{L}_f(\theta, \mathcal{T}_f) &= \frac{1}{T_f}\sum_{\mathbf{x}\in\mathcal{T}_f}\bigg\lVert f\bigg(\mathbf{x},\frac{\partial u^\theta}{\partial x_1},\ldots,\frac{\partial u^\theta}{\partial x_d},\frac{\partial^2 u^\theta}{\partial x_1\partial x_1},\ldots,\frac{\partial^2 u^\theta}{\partial x_1\partial x_d},\ldots,\boldsymbol{\lambda} \bigg) \bigg\rVert_2^2\label{PDE-loss}.
\end{align}
In this form, $f^\theta$ is not directly derived and rather implicitly trained through the PDE loss term, in which the derivatives of $u^{\theta}$ are computed by automatic differentiation. While this framework is quite robust for PDEs on finite spatio-temporal domains, we find that the approach is not well-suited to our task of solving the IVP (\ref{eqn:IVP}) on an unbounded time domain. In this case, the boundary loss term collapses to a single initial condition and terminal condition for our training set with equispaced collocation points $\mathbf{x}\in\mathcal{T}_f$. As a result, the PINN may, but not always, interpolate the training data well within that compact time domain, but we find that it typically does not generalize to open-ended forecasting from other initial conditions. In addition, training by minimizing $\mathcal{L}(\theta, \mathcal{T}_b, \mathcal{T}_f)$, (see \cite{kingma_adam_2017,byrd_limited_1995}), is expensive and slow to converge, especially when compared to NVAR's linear least-squares training. However, the idea that the structure of $f$ can be enforced through a secondary physics-informed neural network does generalize well to NVAR.

\subsection{NVAR differentiation}\label{subsec-nvar-diff}

Recall the NVAR training objective function \linebreak given by Equation (\ref{orig_obj}):
\begin{align*}
	\min_\mathbf{W} \quad w_d g_d(\mathbf{W})  + r g_r(\mathbf{W}).
\end{align*}
We propose the addition of a term analagous to (\ref{PDE-loss}) to measure the correspondence with the model's time derivative and the evaluation of the ODE, (\ref{eqn:IVP}), at each training point. To motivate the form of this term, recall from Section \ref{sec-nvar} the approximation of the definite integral (\ref{I-approx}), letting $s\geq 1$:
\begin{align}\label{short-I-form}
    \mathbf{W}\mathbf{h}(\mathbf{y}_{k}^{p,s})\approx\int_{kh}^{(k+1)h} \mathbf{f}(\mathbf{x}(\xi))d\xi.
\end{align}
By taking the time derivative of each side, we obtain the condition,
\begin{align}
    \frac{d}{dt}\mathbf{W}\mathbf{h}(\mathbf{y}_{k}^{p,s})&\approx\frac{d}{dt}\int_{t+h}^{t} \mathbf{f}(\mathbf{x}(\xi))d\xi \bigg|_{t=kh} \\
    &=\Big(\mathbf{f}\big(\mathbf{x}(t+h)\big)-\mathbf{f}\big(\mathbf{x}(t)\big)\Big) \Big|_{t=kh}\\
    &=\mathbf{f}(\mathbf{x}_{k+1})-\mathbf{f}(\mathbf{x}_{k}).
\end{align}
This indicates that for the training points indexed by $k=a,\ldots,a+T-1$, we should enforce, 
\begin{align}\label{f-approx}
    \frac{d}{dt}\big(\mathbf{W}\mathbf{h}(\mathbf{y}_{k}^{p,s})\big) \approx \mathbf{f}(\mathbf{x}(t_{k+1}))-\mathbf{f}(\mathbf{x}(t_{k})).
\end{align}
The simple expression of NVAR and linear dependence on the trained parameters, $\mathbf{W}$, facilitates the explicit derivation of $\frac{d}{dt}\big(\mathbf{W}\mathbf{h}(\mathbf{y}_{k}^{p,s})\big)$. By the chain rule, we obtain
\begin{align*}
    \frac{d}{dt}\big(\mathbf{W}\mathbf{h}(\mathbf{y}_k^{p,s})\big) &= \mathbf{W}\frac{d}{dt}\mathbf{h}(\mathbf{y}_k^{p,s}) = \mathbf{W}\nabla \mathbf{h}(\mathbf{y}_k^{p,s})\frac{d}{dt}\mathbf{y}_k^{p,s} \\
    &= \mathbf{W}\nabla \mathbf{h}(\mathbf{y}_k^{p,s}) \begin{bmatrix}
        \frac{d}{dt}\mathbf{x}_{k} \\ \frac{d}{dt}\mathbf{x}_{k-s} \\ \vdots \\ \frac{d}{dt}\mathbf{x}_{k-(p-1)s}
    \end{bmatrix}\\
    &= \mathbf{W}\underbrace{\begin{bmatrix}
        \nabla h_1(\mathbf{y}_k^{p,s})^T \\
        \vdots \\
        \nabla h_m(\mathbf{y}_k^{p,s})^T
    \end{bmatrix}}_{\nabla \mathbf{h}(\mathbf{y}_k^{p,s})}\underbrace{\begin{bmatrix}
        \mathbf{f}(\mathbf{x}_k) \\
        \vdots \\
        \mathbf{f}(\mathbf{x}_{k-(p-1)s})
    \end{bmatrix}}_{\mathbf{F}(\mathbf{y}_k^{p,s})}.
\end{align*}
Notice that only $\nabla\mathbf{h}(\mathbf{y}_k^{p,s})$ depends on the model structure, that is, the choice of $\mathbf{h}$, and can therefore be supplied as a function at model initialization. Having obtained $\frac{d}{dt}\big(\mathbf{W}\mathbf{h}(\mathbf{y}_k^{p,s})\big)=\mathbf{W}\nabla \mathbf{h}(\mathbf{y}_k^{p,s})\mathbf{F}(\mathbf{y}_k^{p,s})$, we extend the approximation given by (\ref{f-approx}):
\begin{equation}\label{f-nvar-approx}
\mathbf{f}(\mathbf{x}(t_{k+1}))-\mathbf{f}(\mathbf{x}(t_k))\approx \frac{d}{dt}\big(\mathbf{W}\mathbf{h}(\mathbf{y}_k^{p,s})\big) = \mathbf{W}\nabla \mathbf{h}(\mathbf{y}_k^{p,s})\mathbf{F}(\mathbf{y}_k^{p,s}).
\end{equation}
As a result, we obtain a new update formula for $\mathbf{f}$,
\begin{equation}
\mathbf{f}(\mathbf{x}_{k+1}) = \mathbf{f}(\mathbf{x}_k)+\mathbf{W}\nabla \mathbf{h}(\mathbf{y}_k^{p,s})\mathbf{F}(\mathbf{y}_k^{p,s}),
\label{ode-constraint}
\end{equation}
which provides an additional condition that $\mathbf{W}$ should, mediated by an intermediate map $\nabla\mathbf{h}(\mathbf{y}_k^{p,s})$, transform $(\mathbf{f}(\mathbf{x}_k),\ldots,\mathbf{f}(\mathbf{x}_{k-(p-1)s}))$ to $\mathbf{f}(\mathbf{x}_{k+1})-\mathbf{f}(\mathbf{x}_{k})$. Crucially, the new update formula completely shares the trained parameters, $\mathbf{W}$, with the update for $\mathbf{x}$,
\begin{align}
    \mathbf{x}_{k+1}=\mathbf{x}_k+\mathbf{W}\mathbf{h}(\mathbf{y}_k^{p,s}).
\end{align}
We regard $\mathbf{h}(\mathbf{y}_k^{p,s})\mathbf{F}(\mathbf{y}_k^{p,s})$ as a new state function with prescribed relation to $\mathbf{h}$ and dependence on $\mathbf{f}$. In the context of \cite{raissi_physics-informed_2019}, this is the physics-informed NVAR.

\subsection{Physics-informed training}\label{subsec-physics-informed-training}

To incorporate the condition (\ref{ode-constraint}) into \linebreak NVAR training, we define a new ODE-fit term:
\begin{align}
    &g_o(\mathbf{W}):=\Big\lVert \mathbf{W}\frac{d}{dt}\mathbf{H}-\frac{d}{dt}\mathbf{Z}_{tr} \Big\rVert_F^2,
\end{align}
where
\begin{align}
     &\frac{d}{dt}\mathbf{H}_{tr} = \begin{bmatrix} \nabla \mathbf{h}(\mathbf{y}_{a}^{p,s})\mathbf{F}(\mathbf{y}_{a}^{p,s}) & \ldots & \nabla \mathbf{h}(\mathbf{y}_{a+T-1}^{p,s})\mathbf{F}(\mathbf{y}_{a+T-1}^{p,s})\end{bmatrix},\\
     &\frac{d}{dt}\mathbf{Z}_{tr} = \begin{bmatrix}
         \frac{d}{dt}\mathbf{z}_{a} & \ldots & \frac{d}{dt}\mathbf{z}_{a+T-1}
     \end{bmatrix}.
\end{align}
Recall from Equation \cref{z-nvar-tgt} that for $k=0,\ldots,T-1$, $\mathbf{z}(t_{a+k})=\mathbf{x}(t_{a+k+1})-\mathbf{x}(t_{a+k})$, so $\frac{d}{dt}\mathbf{z}(t_{a+k})=\mathbf{f}\big(\mathbf{x}(t_{a+k+1})\big)-\mathbf{f}\big(\mathbf{x}(t_{a+k})\big)$. The new training objective is then given by 
\begin{align}\label{ode_obj}
	\min_{\mathbf{W}\in\mathbb{R}^{d\times m}} \quad \underbrace{w_d g_d(\mathbf{W})}_{\text{data-fit}} + \underbrace{w_o g_o(\mathbf{W})}_{\text{ODE-fit}} + \underbrace{r g_r(\mathbf{W})}_{\text{regularization}}.
\end{align}
Similar to (\ref{orig_obj}) and (\ref{equiv_obj}), we formulate an equivalent objective,
\begin{align}\label{equiv_ode_obj}
	\min_{\mathbf{W}\in\mathbb{R}^{d\times m}} \quad \big\lVert \mathbf{L}(\mathbf{W})\big\rVert_F^2,
\end{align}
where
\begin{align}\label{equiv_ode_L}
	\mathbf{L}(\mathbf{W}) =  \mathbf{W} \begin{bmatrix}
	    \sqrt{w_d}\ \mathbf{H} \\
	    \sqrt{w_o}\ \frac{d}{dt}\mathbf{H} \\
        \sqrt{r}\ \mathbf{I}
	\end{bmatrix} - \begin{bmatrix}
	    \sqrt{w_d}\ \mathbf{Z}_{tr} \\
	    \sqrt{w_o}\ \frac{d}{dt}\mathbf{Z}_{tr} \\
        \mathbf{0}
	\end{bmatrix}.
 \end{align}
 The objective function remains a linear least-squares problem in $\mathbf{W}$ and no adjustment to the solution process is necessary.

\subsection{Illustrative example}\label{subsec-example}

We first give an example of a simple state function, $\mathbf{h}$, in particular the form resulting from the derivation from ESNs in \cite{bollt_explaining_2021}. Let
\begin{align}\label{eqn:h}
    \mathbf{h}(\mathbf{y}_k^{p,s})&=\begin{bmatrix}
        1 \\
        \mathbf{y}_k^{p,s} \\
        \mathbf{h}_{nonlin}(\mathbf{y}_k^{p,s})
    \end{bmatrix},
\end{align}
in which $\mathbf{h}$ is the concatenation of a constant bias entry, a linear portion , $\mathbf{y}_k^{p,s}\in\mathbb{R}^{pd}$, and a nonlinear portion, $\mathbf{h}_{nonlin}(\mathbf{y}_k^{p,s})$, where $\mathbf{h}_{nonlin}:\mathbb{R}^{pd}\mapsto \mathbb{R}^{n}$. Here, $n$ depends on the choice of nonlinearity. As a result, the total dimension of the state vector is $m=1+pd+n$. For example, a common choice of nonlinear function is all possible quadratic monomials of the linear terms without duplicates. Borrowing notation from \cite{bollt_explaining_2021}, we define
\begin{align}
\begin{split}
    \mathbf{p}_2:\mathbb{R}^n\times\mathbb{R}^n &\rightarrow \mathbb{R}^{n(n+1)/2} \\
    \mathbf{p}_2(\mathbf{v},\mathbf{w}) &\mapsto [v_1w_1 \ \ldots \ v_1w_n \ v_2w_1 \ \ldots \ v_2w_n \ \ldots v_nw_n
    ]^T,
\end{split}
\end{align}
and choose
\begin{align}\label{h_nonlin}
	\mathbf{h}_{nonlin}(\mathbf{y}_k^{p,s})&= \mathbf{p}_2(\mathbf{y}_k^{p,s},\mathbf{y}_k^{p,s}) =\begin{bmatrix}
		x_{k,1}^2 \\
        x_{k,1}x_{k,2} \\
		\vdots \\
		x_{k-(p-1),d-1}x_{k-(p-1),d} \\
		x_{k-(p-1),d}^2
	\end{bmatrix}.
\end{align}
In this case, $n=\frac{1}{2}(pd(pd+1))$. We also demonstrate the form of $\nabla \mathbf{h}$ in this minimal example. Let $\mathbf{x}=(x_1,x_2)$ and
\begin{align}
    \mathbf{h}(\mathbf{x}) &= \begin{bmatrix}
        1 \\
        \mathbf{x} \\
        \mathbf{p}_2(\mathbf{x},\mathbf{x})
    \end{bmatrix} = \begin{bmatrix}
        1 \\ x_1 \\ x_2 \\ x_1^2 \\ x_1 x_2 \\ x_2^2
    \end{bmatrix},
\end{align}
then,
\begin{align}\label{grad-h-example}
    \nabla\mathbf{h}(\mathbf{x}) = \begin{bmatrix}
        \frac{\partial}{\partial x_1}1 & \frac{\partial}{\partial x_2}1 \\
        \frac{\partial}{\partial x_1}x_1 & \frac{\partial}{\partial x_2}x_1 \\
        \frac{\partial}{\partial x_1}x_2 & \frac{\partial}{\partial x_2}x_2 \\
        \frac{\partial}{\partial x_1}x_1^2 & \frac{\partial}{\partial x_2}x_1^2 \\
        \frac{\partial}{\partial x_1}x_1 x_2 & \frac{\partial}{\partial x_2}x_1 x_2 \\
        \frac{\partial}{\partial x_1}x_2^2 & \frac{\partial}{\partial x_2}x_2^2 \\
    \end{bmatrix} = \begin{bmatrix}
        0 & 0 \\
        1 & 0 \\
        0 & 1 \\
        2x_1 & 0 \\
        x_2 & x_1 \\
        0 & 2x_2 \\
    \end{bmatrix}.
\end{align}

\section{Numerical experiments}\label{sec-num-exp}

In this Section, we describe numerical experiments to evaluate the performance of piNVAR to predict various dynamical systems. Subsection \ref{subsec-test-problems} provides a description of the test problems. The cross-validation testing routine is detailed in Subsection \ref{subsec-cross-val}. Subsection \ref{subsec-eval-metrics} reviews the two evaluation metrics and Subsection \ref{subsec-model-params} gives the model parameters. Julia 1.8.0 \cite{Julia-2017} was used for all implementation.

\subsection{Test problems}\label{subsec-test-problems}

We use three ordinary differential equations. They encompass linear, nonlinear, and chaotic dynamics both with known and unknown analytical solutions.

\subsubsection{Undamped spring}

The undamped spring equation is ubiquitous in basic ODE theory. It is a linear second order equation in one dimension given by
\begin{align}\label{second-order-spring}
    x''(t)+kx(t)&=0,\\
    x(0)=0, \\
    x'(0)=1,
\end{align}
where $k$ is a spring constant. For the numerical tests presented here, we take $k=3$. The spring equation has the simple exact solution,
\begin{align}
    x(t)=\sin(\sqrt{k}t).
\end{align}
Since we only consider first-order ODEs in our formulation of physics-informed NVAR, we convert the above second-order one-dimensional equation into a first-order system in two variables:
\begin{align}
    \mathbf{x}(t) &= \begin{bmatrix}
        x_1(t) \\
        x_2(t)
    \end{bmatrix} = \begin{bmatrix}
        x_1(t) \\
        x_1'(t)
    \end{bmatrix}, \\
    \mathbf{x}'(t) &= \begin{bmatrix}
        x_1'(t)\\
        x_2'(t)
    \end{bmatrix} = \begin{bmatrix}
        x_2(t) \\
        -kx_1(t)
    \end{bmatrix}.
\end{align}
The exact solution for the first-order system is
\begin{align}
    \mathbf{x}(t) = \begin{bmatrix}
			\sin(\sqrt{k}t) \\
			\sqrt{k}\cos(\sqrt{k}t)
		\end{bmatrix}.
\end{align}

\subsubsection{Lotka-Volterra}

The Lotka-Volterra equation is a nonlinear first order system in two variables that describes the population dynamics of a prey species, $x_1$, and predator species, $x_2$. Its parameters describe the inter-species interactions: $a$ is the per-capita prey growth rate, $b$ is the per-capita prey death rate from predation, $c$ is the per-capita predator growth rate from predation, and $d$ is the per-capita predator death rate. The system is given by
\begin{align}
    \mathbf{x}'(t)&=\begin{bmatrix}
        x_1'(t) \\
        x_2'(t)
    \end{bmatrix}=\begin{bmatrix}
        a x_1(t)-b x_1(t)x_2(t) \\
        c x_1(t)x_2(t)-d x_2(t)
    \end{bmatrix}, \\
    \mathbf{x}(0) &= \begin{bmatrix}
        1.0 \\
        0.25
    \end{bmatrix}.
\end{align}
For the tests presented here, we choose $a=0.25$, $b=1.0$, $c=0.5$, and $d=0.125$. We generate an approximate numerical solution using the RK4 scheme using time step $h=\num{1e-5}$ and down-sampling to an effective time-step of $h=\num{1e-1}$ by selecting every thousandth point.

\subsubsection{Lorenz}
The nonlinear and chaotic Lorenz system of ODEs models atmospheric convection \cite{lorenz_deterministic_1963}. We constrain ourselves to canonical parameter values that yield chaotic behavior,
\begin{align}\label{eqn:lorenz}
    \mathbf{x}'(t)=\begin{bmatrix}
        x_1'(t) \\ x_2'(t) \\ x_3'(t)
    \end{bmatrix}=\begin{bmatrix}
        10(x_2-x_1) \\ x_1(28-x_3)-x_2 \\ x_1x_2-\frac{8}{3}x_3
    \end{bmatrix},\quad \mathbf{x}(0)=\begin{bmatrix}
        -3 \\-3 \\ 28
    \end{bmatrix}.
\end{align}
We generate an approximate numerical solution using RK4 using time step $h=\num{1e-5}$ and down-sampling to an effective time-step of $h=\num{1e-3}$ by selecting every hundredth point. It is important to note that referring to this numerical solution as trustworthy reference data is inherently flawed and does not necessarily represent the exact trajectory as $t\rightarrow \infty$ given $\mathbf{x}(0)$ due to the chaotic nature of the ODE. However, it does serve as a fixed training set and data-driven benchmark for different model parameters.

\subsection{Cross-validation routine}\label{subsec-cross-val}

We are interested in evaluating the ability of an NVAR model to produce good predictions for multiple initial conditions, given only a single small training set. In other words, we fix a training set, apply the trained model to a handful of differing initial conditions, and measure the minimum, median, or maximum results of a particular evaluation metric across those trials. For each test problem, we generate 100,000 reference data points by the procedures outlined in the preceding subsections. We partition each data set into training and testing sets using index sets which are the same for all problems. The training indices, $\mathcal{I}_{tr}$, are the 1,500 points indexed by $k=2{,}001,\ \ldots,\ 3{,}500$. We set five disjoint intervals, each 10,000 points which are disjoint from the training set, for recursive prediction and cross-validation. 
\begin{align}
\begin{split}
    \mathcal{I}_{test1} &= \{10{,}001,\ \ldots,\ 20{,}000\},\\
    \mathcal{I}_{test2} &= \{20{,}001,\ \ldots,\ 30{,}000\},\\
    \mathcal{I}_{test3} &= \{30{,}001,\ \ldots,\ 40{,}000\},\\
    \mathcal{I}_{test4} &= \{40{,}001,\ \ldots,\ 50{,}000\},\\
    \mathcal{I}_{test5} &= \{50{,}001,\ \ldots,\ 60{,}000\}.
\end{split}
\end{align}

\subsection{Evaluation metrics}\label{subsec-eval-metrics}

\subsubsection{Valid time}
\begin{figure}[h]
    \centering
    \includegraphics[scale = 0.60]{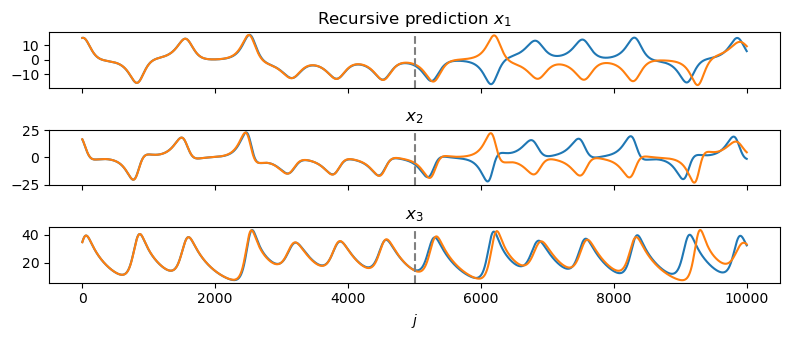}
    \caption[Valid time for example recursive prediction]{Recursive prediction of system state $\mathbf{x}=[x_1\ x_2\ x_3]^\top$ for the Lorenz system. Blue: model prediction, orange: target, dashed line: valid time for $M=0.01$. Here, $t_{valid}=5,005$}
    \label{pred_fig}
\end{figure}

For the chaotic Lorenz ODE, we find that a recursive prediction often matches the reference data well until accumulated error causes the predicted trajectory to abruptly diverge. As a result, popular metrics like root mean-squared error (RMSE) do not effectively capture the length of predictive validity of a model. Instead, we adopt the notion of valid time from \cite{rodrigues_physics-informed_2019} as our primary data-driven evaluation metric. Valid time measures the number of time steps until the relative error of the recursively predicted system state exceeds a specified threshold, $M$. Let the sequence of predicted points be $\mathbf{y}_j$ and reference data $\mathbf{x}_j$, $1\leq j \leq T_{test}$.  We define the valid time:
\begin{align}
    t_{valid}:=\min_j \bigg\{j\ s.t.\ \frac{\lVert \mathbf{y}_j-\mathbf{x}_j\rVert_2^2}{\lVert \mathbf{x}_j\rVert_2^2}\geq M \bigg\},
\end{align}
Figure \ref{pred_fig} gives an example of recursive prediction, testing, and valid time for an example model with $M=0.01$. In all of the numerical experiments presented here, we choose $M=\num{1e-4}$.

\subsubsection{Discrete energy functional}
While valid time measures the ability of an NVAR model to generate an approximate solution which closely matches desired reference data, we propose a physics-focused metric to measure the quality of the approximate physics-informed NVAR solution in terms of adherence to the underlying ODE. This is particularly important when the reference data differs from the exact IVP solution, which is most relevant in the case of the chaotic Lorenz system. Borrowing from error analysis of physics-informed neural network solutions to PDEs \cite{zeinhofer_unified_2024}, we formulate a discrete approximation of an energy functional over a prescribed prediction interval, $[0,t_T]$. Recall that each ODE is described by
\begin{align}
    \frac{d}{dt}\mathbf{x}(t) &= \mathbf{f}(\mathbf{x}(t)) \iff \mathbf{f}(\mathbf{x}(t))-\frac{d}{dt}\mathbf{x}(t) = 0.
\end{align}
We define an energy functional,
\begin{align}
    E(T) &:= \frac{1}{2}\int_{0}^{t_T} \Big( \mathbf{f}(\mathbf{x}(t))-\frac{d}{dt}\mathbf{x}(t) \Big)^2 dt,
\end{align}
and approximate it using simple right-endpoint rectangular quadrature by
\begin{align}
    E(T) &\approx \frac{1}{2} \sum_{k=0}^{T-1} h\Big( \mathbf{f}(\mathbf{x}(t_{k+1}))-\frac{d}{dt}\mathbf{x}(t_{k+1})\Big)^2.
\end{align}
Finally, we recall from \eqref{f-nvar-approx} that the physics-informed NVAR approximately enforces the relationship,
\begin{align}
    \frac{d}{dt}\mathbf{x}(t_{k+1}) &= \mathbf{f}(\mathbf{x}(t_{k+1})) \approx \mathbf{f}(\mathbf{x}_{k+1}) \approx \mathbf{f}(\mathbf{x}_{k})+\mathbf{W}\nabla\mathbf{h}(\mathbf{y}_k^{p,s})\mathbf{F}(\mathbf{y}_k^{p,s}),
\end{align}
which we substitute to obtain
\begin{align}
    E(T) &\approx \frac{h}{2} \sum_{k=0}^{T-1} \Big( \mathbf{f}(\mathbf{x}_{k+1})-\big(\mathbf{f}(\mathbf{x}_{k})+\mathbf{W}\nabla\mathbf{h}(\mathbf{y}_k^{p,s})\mathbf{F}(\mathbf{y}_k^{p,s})
    \big)\Big)^2 =: E_h(T).
\end{align}
Therefore, $E_h(T)$ measures the approximate least-squares deviation from satisfying the governing ODE, $\mathbf{x}'(t)-\mathbf{f}(\mathbf{x}(t))=\mathbf{0}$, over the prediction interval $[0,t_T]$. Crucially, this metric is independent of any reference data.

\subsection{Model parameters}\label{subsec-model-params}

We test the effectiveness of piNVAR for three state functions, each of which is motivated by low-degree polynomial approximation. Define the space of multivariate polynomials with $n$ total maximal degree by,
\begin{align}
    \mathcal{P}_n(\mathbf{x}) := \bigg\{\prod_{i=1}^d x_i^{\alpha_i}\ :\ \sum_{i=1}^d|\alpha_i|\leq n\ \bigg\}
\end{align}
By definition, $\mathcal{P}_n(\mathbf{x})$ is the monomial basis for itself. Denote the $n$-degree Chebyshev polynomial of the first kind by $T_n(x)$. Then, the Chebyshev basis for $\mathcal{P}_n(\mathbf{x})$ is
\begin{align}
    \bigg\{\prod_{i=1}^d T_{\alpha_i}(x_i)\ :\ \sum_{i=1}^d|\alpha_i|\leq n\ \bigg\}.
\end{align}
The three state functions are described in Table \ref{state-funcs}.
\begin{table}[H]
\begin{center}
\renewcommand{\arraystretch}{0.75}
\begin{tabular}{|c|l|}
\hline
Name & Description \\
\hline
$\mathbf{h}_1$ & Monomial basis for $\mathcal{P}_2(\mathbf{x})$ with smooth support \\
$\mathbf{h}_2$ & Monomial basis for $\mathcal{P}_2(\mathbf{x})$ with non-smooth support \\
$\mathbf{h}_3$ & Chebyshev basis for $\mathcal{P}_2(\mathbf{x})$ with non-smooth input transformation \\
\hline
\end{tabular}
\caption[Tested state functions.]{Tested state functions.}\label{state-funcs}\end{center}\end{table}

The first two functions, $\mathbf{h}_1$ and $\mathbf{h}_2$, are as seen in Subsection \ref{subsec-example}, each with an adjustment to add local support to the otherwise unbounded monomials. The model is trained on sampled data points from the measured historical trajectory of the ODE, so these coordinates are effectively bounded within the feasible domain of the system. However, when the model is used for recursive prediction, the forecasted trajectory can drift outside the true feasible domain and may blow up. One potential remedy, easily implemented, is to define a local support coefficient function, $\phi:\mathbf{R}\mapsto[0,1]$, which rapidly goes to zero outside a desired compact domain. When multiplied into the entries of $\mathbf{h}$, the local support coefficients reduce model inputs from coordinates which stray outside the desired attractor domain. 

The general form of $\mathbf{h}_1$ and $\mathbf{h}_2$ is given below. The only difference between the two is the particular choice of $\phi$.
\begin{align}\label{test-h}
    \mathbf{h}(\mathbf{y}_k^{p,s})&=\begin{bmatrix}
        1 \\
        \boldsymbol{\phi}_y(\mathbf{y}_k^{p,s})\odot \mathbf{y}_k^{p,s} \\
        \mathbf{h}_{nonlin}\big(\boldsymbol{\phi}_y(\mathbf{y}_k^{p,s})\odot \mathbf{y}_k^{p,s}\big)
    \end{bmatrix}, \\
    \mathbf{h}_{nonlin}(\mathbf{y})&= \mathbf{p}_2(\mathbf{y},\mathbf{y} ),
\end{align}
where
\begin{align}
    \boldsymbol{\phi}_y(\mathbf{y}_k^{p,s}) &= \begin{bmatrix}
        \boldsymbol{\phi}(\mathbf{x}_k) \\
        \vdots \\
        \boldsymbol{\phi}(\mathbf{x}_{k-(p-1)s})
    \end{bmatrix} = \begin{bmatrix}
        \phi_1(x_{k,1}) \\
        \vdots \\
        \phi_d(x_{k,d}) \\
        \vdots \\
        \phi_1(x_{k-(p-1)s,1}) \\
        \vdots \\
        \phi_d(x_{k-(p-1)s,d})
    \end{bmatrix}.
\end{align}
Here, $\phi_i$ denotes the support coefficient function corresponding to the $i^{th}$ coordinate. It is reasonable to expect that each coordinate, $x_i$, $i=1,\ldots,d$, occupies a different desired domain, so we implement support functions with coordinate-specific parameters, $\phi_i(\cdot) := \phi(\cdot\ ;\boldsymbol{\lambda}_i)$.

For $\mathbf{h}_1$, $\phi$ is based on a hyperbolic tangent that smoothly transitions from 1 to 0 across a specified radius $r\in\mathbb{R}$ from the center, $t\in\mathbb{R}$, with specified sharpness $\xi \in\mathbb{Z}^+$,
\begin{align}\label{eqn:bump}
    \phi(x;\boldsymbol{\lambda}) = \phi(x;\xi, r, t):=\frac{1}{2}\Big(1+\tanh\big((\xi+\frac{1}{4})\pi(r^2-(x-t)^2)\big)\Big).
\end{align}
For each coordinate, the center is $t=0$, the sharpness is $\xi=5$, and each radius $r_i$ is set to 110\% of the maximum coordinate-wise magnitude in the reference data,
\begin{align}
    r_i = 1.1\max_{1\leq k \leq T} |x_{k,i}|.
\end{align}
As a result, $\boldsymbol{\lambda}_i = (5, r_i, 0)$.

For $\mathbf{h}_2$, the smooth support coefficient is replaced with a piecewise linear function with transitions over specified intervals, $[a,b]$ and $[c,d]$.
\begin{align}\label{eqn:nonsmooth-bump}
    \phi(x;\boldsymbol{\lambda}) = \phi(x;a,b,c,d):= \begin{cases}
        0 \quad & x \leq a, \\
        \frac{x-a}{b-a} \quad & a < x \leq b, \\
        1 \quad & b < x \leq c, \\
        1-\frac{x-c}{d-c} \quad & c < x \leq d, \\
        0 \quad & x > d.
    \end{cases}
\end{align}
Here, we choose the coordinate-specific parameters $a_i=-r_i$, $b_i=-0.95r_i$, $c_i=0.95r_i$, and $d_i=r_i$, and thus $\boldsymbol{\lambda}_i = (-r_i, -0.95r_i, 0.95r_i, r_i)$.

A minimal example is illustrative of the form of $\nabla \mathbf{h}$ for $\mathbf{h}_1$ and $\mathbf{h}_2$. Let 
\begin{align}
    \mathbf{h}=\begin{bmatrix}
        1  & x_1\phi_1(x_1) & x_2\phi_2(x_2) & x_1^2\phi_1(x_1)^2 & x_1 x_2 \phi_1(x_1)\phi_2(x_2) & x_2^2\phi_2(x_2)^2
    \end{bmatrix}^\top.
\end{align}
Then,
\begin{align}\label{supp-grad-h-example}
    \nabla&\mathbf{h} = \begin{bmatrix}
        0 & 0 \\
        \phi_1(x_1)+x_1\phi_1(x_1)' & 0 \\
        0 & \phi_2(x_2)+x_2\phi_2(x_2)' \\
        2x_1\phi_1(x_1)^2+2x_1^2\phi_1(x_1)\phi_1(x_1)' & 0 \\
        x_2\phi_1(x_1)\phi_2(x_2)+x_1x_2\phi_1(x_1)'\phi_2(x_2) & x_1\phi_1(x_1)\phi_2(x_2)+x_1x_2\phi_1(x_1)\phi_2(x_2)' \\
        0 & 2x_2\phi_2(x_2)^2+2x_2^2\phi_2(x_2)\phi_2(x_2)' \\
    \end{bmatrix}
\end{align}
Plots of both smooth and non-smooth support coefficient functions are included below.
\begin{figure}[H]
\centering
\begin{minipage}{0.55\textwidth}
  \centering
\includegraphics[height=2.15in]{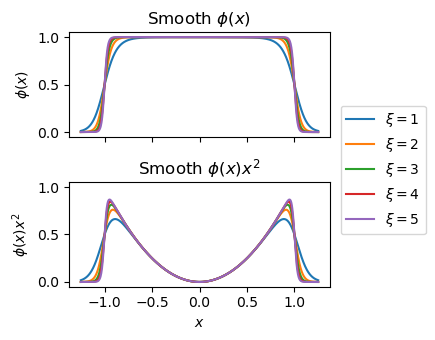}

\end{minipage}%
\begin{minipage}{0.45\textwidth}
  \centering
\includegraphics[height=2.15in]{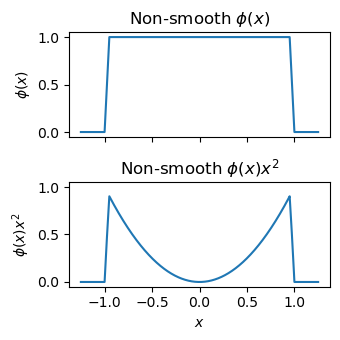}
\end{minipage}%

\caption[Support functions]{Support coefficient functions. Left: smooth support for $m=1.0$, $t=0.0$, and sharpness $\xi=1,\ldots,5$, right: non-smooth support for $a=-1$, $b=-0.85$, $c=0.85$, $d=1.0$} \label{support-coeffs}
\end{figure}

While $\mathbf{h}_1$ and $\mathbf{h}_2$ ensure the boundedness of the recursive NVAR prediction by preventing unbounded functions from taking input from out-of-bounds coordinates, $\mathbf{h}_3$ achieves boundedness naturally through Chebyshev polynomials of the first kind, $T_n(x)=\cos(n\arccos(x))$ for $x\in[-1,1]$. Some input transformation is necessary to ensure $T_n$ is defined for all likely predicted points. For $a<b$, we define,
\begin{align}
    \phi_{ab}(x) := \begin{cases}
        \frac{a}{x} \quad & x < a, \\
        1 \quad & a \leq x < b, \\
        \frac{b}{x} \quad & x \geq b, \\
    \end{cases}
\end{align}
which yields the piecewise linear composite mapping,
\begin{align}
    \phi_{ab}(x)x = \begin{cases}
        a \quad & x < a, \\
        x \quad & a \leq x < b, \\
        b \quad & x \geq b. \\
    \end{cases}
\end{align}
Then, we define the affine linear transformation,
\begin{align}
    A_{ab}(x):= \frac{a+b-2x}{a-b},
\end{align}
which, composed with the previous mapping gives,
\begin{align}
    \Lambda_{ab}(x) := A_{ab}\big(\phi_{ab}(x)x \big) = \begin{cases}
        -1 \quad & x < a, \\
        \frac{a+b-2x}{a-b} \quad & a \leq x < b, \\
        1 \quad & x \geq b. \\        
    \end{cases}
\end{align}
Now, $T_{n,ab}(x):=T_n(\Lambda_{ab}(x))$ is defined for all real-valued $x$. To promote clarity, we replace the references to the coordinate-wise parameters $a_i$ and $b_i$ and use the notation $T_{n,i}(x_i)$ in place of $T_{n,a_ib_i}(x_i)$. To connect this back to the form of $\mathbf{h}_3$, let us suppose $\mathbf{x}=(x_1,x_2)$. The multivariate Chebyshev basis for $\mathcal{P}_2(x_1,x_2)$ is
\begin{align}
    \begin{Bmatrix}
        T_0(x_1)T_0(x_2),\\
        T_0(x_1)T_1(x_2),\\
        T_1(x_1)T_0(x_2),\\
        T_2(x_1)T_0(x_2),\\
        T_1(x_1)T_1(x_2),\\
        T_0(x_1)T_2(x_2)
    \end{Bmatrix}.
\end{align}
In each term, we replace $x_i$ with $\Lambda_{a_ib_i}(x_i)$ and define
\begin{align}
    \mathbf{h}_3(\mathbf{x}) &= \mathbf{h}(x_1,x_2) = \begin{bmatrix}
        T_{0,1}(x_1)T_{0,2}(x_2)\\
        T_{0,1}(x_1)T_{1,2}(x_2)\\
        T_{1,1}(x_1)T_{0,2}(x_2)\\
        T_{2,1}(x_1)T_{0,2}(x_2)\\
        T_{1,1}(x_1)T_{1,2}(x_2)\\
        T_{0,1}(x_1)T_{2,2}(x_2)
    \end{bmatrix},
\end{align}

The NVAR input parameters are $p=10$ and $s=1$. The resulting model size, e.g. the number of entries of $\mathbf{W}$, is 250 for $\mathbf{h}_1$ and $\mathbf{h}_2$ and 231 for $\mathbf{h}_3$ for test problems with $d=2$ and 525 for $\mathbf{h}_1$ and $\mathbf{h}_2$ and 496 for $\mathbf{h}_3$ for test problems with $d=3$.

\subsection{Results}\label{subsec-results}

To evaluate the impact of physics-informed training, on a range of models, we vary regularization $r$ and the ODE training weight $w_o$ across a discrete dictionary of values and assess the median valid time, taking $M=\num{1e-4}$, and median discrete energy, $E_h(10,000)$, across the five separate trials for each model. Tables \ref{spring-vt-table}-\ref{lorenz-vt-table} display the median valid time results for the spring, Lotka-Volterra, and Lorenz test problems, respectively. Table \ref{spring-vt-table} shows that most model configurations achieve a median valid time of 10,000 time steps, the maximum possible in our test. Despite the broadly strong performance, we do note a tradeoff between regularization and median valid time, with larger values of $r$ producing lower valid times. However, meaningful physics-informed training, measured by $w_o$, provides a strong mitigant. Table \ref{lv-vt-table} shows generally lower valid times and an even stronger inverse relationship between regularization and valid time, suggesting that the Lotka-Volterra test problem is more difficult for NVAR. We also see that the mitigation effect of $w_o$ is substantially weaker. This result is interesting when contrasted with the Lorenz test problem, which, due to its chaotic nature, is typically considered more difficult to address with numerical methods. With the exception of the poorly regularized models, Table \ref{lorenz-vt-table} shows generally greater valid times, weaker deterioration from regularization, and strong performance improvement from physics-informed training. In all but three cases, the superior model for each state function, $\mathbf{h}_i$, and regularization value is achieved by either 0.5 or 1.0 ODE training weight.

We also evaluate the cross-validation results using the discrete energy functional, $E_h(10,000)$. The median value of $E_h(10,000)$ across the five 10,000 time step prediction intervals is shown in Tables \ref{spring-et-table}-\ref{lorenz-et-table} for the spring, Lotka-Volterra, and Lorenz test problems, respectively. Starting with the spring problem, we see a substantial difference between state functions for the first time, with $\mathbf{h}_2$ and $\mathbf{h}_3$ achieving up to 21 orders of magnitude improvement on median $E_h$ compared to $\mathbf{h}_1$ for small regularization. The difference shrinks to just one to three orders of magnitude for $r=\num{1e-1}$. For all state functions, large ODE training weight tends to yield the superior model, with the performance improvement most pronounced as regularization increases. 

\newpage

\begin{table}[H]
\begin{center}
\renewcommand{\arraystretch}{0.75}
\begin{tabular}{|c|c|cccccc|}
\hline  
 & $r$ \ \textbackslash \ $w_o$ & 0 & $\num{1e-4}$ & $\num{1e-2}$ & $\num{1e-1}$ & 0.5 & 1.0\\ 
\hline
\multirow{5}{*}{$\mathbf{h}_1$} & $\num{1e-12}$ & \cellcolor[HTML]{C6EFCE}{\color[HTML]{006100} 10000} & \cellcolor[HTML]{C6EFCE}{\color[HTML]{006100} 10000} & \cellcolor[HTML]{C6EFCE}{\color[HTML]{006100} 10000} & \cellcolor[HTML]{C6EFCE}{\color[HTML]{006100} 10000} & \cellcolor[HTML]{C6EFCE}{\color[HTML]{006100} 10000} & \cellcolor[HTML]{C6EFCE}{\color[HTML]{006100} 10000} \\
 & $\num{1e-8}$ & \cellcolor[HTML]{C6EFCE}{\color[HTML]{006100} 10000} & \cellcolor[HTML]{C6EFCE}{\color[HTML]{006100} 10000} & \cellcolor[HTML]{C6EFCE}{\color[HTML]{006100} 10000} & \cellcolor[HTML]{C6EFCE}{\color[HTML]{006100} 10000} & \cellcolor[HTML]{C6EFCE}{\color[HTML]{006100} 10000} & \cellcolor[HTML]{C6EFCE}{\color[HTML]{006100} 10000} \\
 & $\num{1e-4}$ & \cellcolor[HTML]{C6EFCE}{\color[HTML]{006100} 10000} & \cellcolor[HTML]{C6EFCE}{\color[HTML]{006100} 10000} & \cellcolor[HTML]{C6EFCE}{\color[HTML]{006100} 10000} & \cellcolor[HTML]{C6EFCE}{\color[HTML]{006100} 10000} & \cellcolor[HTML]{C6EFCE}{\color[HTML]{006100} 10000} & \cellcolor[HTML]{C6EFCE}{\color[HTML]{006100} 10000} \\
 & $\num{1e-2}$ & 1790 & 1797 & 5228 & \cellcolor[HTML]{C6EFCE}{\color[HTML]{006100} 10000} & \cellcolor[HTML]{C6EFCE}{\color[HTML]{006100} 10000} & \cellcolor[HTML]{C6EFCE}{\color[HTML]{006100} 10000} \\
 & $\num{1e-1}$ & 1309 & 1312 & 1522 & 5234 & \cellcolor[HTML]{C6EFCE}{\color[HTML]{006100} 10000} & \cellcolor[HTML]{C6EFCE}{\color[HTML]{006100} 10000} \\
 \hline
\multirow{5}{*}{$\mathbf{h}_2$} & $\num{1e-12}$ & \cellcolor[HTML]{C6EFCE}{\color[HTML]{006100} 10000} & \cellcolor[HTML]{C6EFCE}{\color[HTML]{006100} 10000} & \cellcolor[HTML]{C6EFCE}{\color[HTML]{006100} 10000} & \cellcolor[HTML]{C6EFCE}{\color[HTML]{006100} 10000} & \cellcolor[HTML]{C6EFCE}{\color[HTML]{006100} 10000} & \cellcolor[HTML]{C6EFCE}{\color[HTML]{006100} 10000} \\
 & $\num{1e-8}$ & \cellcolor[HTML]{C6EFCE}{\color[HTML]{006100} 10000} & \cellcolor[HTML]{C6EFCE}{\color[HTML]{006100} 10000} & \cellcolor[HTML]{C6EFCE}{\color[HTML]{006100} 10000} & \cellcolor[HTML]{C6EFCE}{\color[HTML]{006100} 10000} & \cellcolor[HTML]{C6EFCE}{\color[HTML]{006100} 10000} & \cellcolor[HTML]{C6EFCE}{\color[HTML]{006100} 10000} \\
 & $\num{1e-4}$ & \cellcolor[HTML]{C6EFCE}{\color[HTML]{006100} 10000} & \cellcolor[HTML]{C6EFCE}{\color[HTML]{006100} 10000} & \cellcolor[HTML]{C6EFCE}{\color[HTML]{006100} 10000} & \cellcolor[HTML]{C6EFCE}{\color[HTML]{006100} 10000} & \cellcolor[HTML]{C6EFCE}{\color[HTML]{006100} 10000} & \cellcolor[HTML]{C6EFCE}{\color[HTML]{006100} 10000} \\
 & $\num{1e-2}$ & \cellcolor[HTML]{C6EFCE}{\color[HTML]{006100} 10000} & \cellcolor[HTML]{C6EFCE}{\color[HTML]{006100} 10000} & \cellcolor[HTML]{C6EFCE}{\color[HTML]{006100} 10000} & \cellcolor[HTML]{C6EFCE}{\color[HTML]{006100} 10000} & \cellcolor[HTML]{C6EFCE}{\color[HTML]{006100} 10000} & \cellcolor[HTML]{C6EFCE}{\color[HTML]{006100} 10000} \\
 & $\num{1e-1}$ & 1409 & 1413 & 3183 & \cellcolor[HTML]{C6EFCE}{\color[HTML]{006100} 10000} & \cellcolor[HTML]{C6EFCE}{\color[HTML]{006100} 10000} & \cellcolor[HTML]{C6EFCE}{\color[HTML]{006100} 10000} \\
 \hline
\multirow{5}{*}{$\mathbf{h}_3$} & $\num{1e-12}$ & \cellcolor[HTML]{C6EFCE}{\color[HTML]{006100} 10000} & \cellcolor[HTML]{C6EFCE}{\color[HTML]{006100} 10000} & \cellcolor[HTML]{C6EFCE}{\color[HTML]{006100} 10000} & \cellcolor[HTML]{C6EFCE}{\color[HTML]{006100} 10000} & \cellcolor[HTML]{C6EFCE}{\color[HTML]{006100} 10000} & \cellcolor[HTML]{C6EFCE}{\color[HTML]{006100} 10000} \\
 & $\num{1e-8}$ & \cellcolor[HTML]{C6EFCE}{\color[HTML]{006100} 10000} & \cellcolor[HTML]{C6EFCE}{\color[HTML]{006100} 10000} & \cellcolor[HTML]{C6EFCE}{\color[HTML]{006100} 10000} & \cellcolor[HTML]{C6EFCE}{\color[HTML]{006100} 10000} & \cellcolor[HTML]{C6EFCE}{\color[HTML]{006100} 10000} & \cellcolor[HTML]{C6EFCE}{\color[HTML]{006100} 10000} \\
 & $\num{1e-4}$ & \cellcolor[HTML]{C6EFCE}{\color[HTML]{006100} 10000} & \cellcolor[HTML]{C6EFCE}{\color[HTML]{006100} 10000} & \cellcolor[HTML]{C6EFCE}{\color[HTML]{006100} 10000} & \cellcolor[HTML]{C6EFCE}{\color[HTML]{006100} 10000} & \cellcolor[HTML]{C6EFCE}{\color[HTML]{006100} 10000} & \cellcolor[HTML]{C6EFCE}{\color[HTML]{006100} 10000} \\
 & $\num{1e-2}$ & 8789 & 8799 & \cellcolor[HTML]{C6EFCE}{\color[HTML]{006100} 10000} & \cellcolor[HTML]{C6EFCE}{\color[HTML]{006100} 10000} & \cellcolor[HTML]{C6EFCE}{\color[HTML]{006100} 10000} & \cellcolor[HTML]{C6EFCE}{\color[HTML]{006100} 10000} \\
 & $\num{1e-1}$ & 1356 & 1359 & 1729 & \cellcolor[HTML]{C6EFCE}{\color[HTML]{006100} 10000} & \cellcolor[HTML]{C6EFCE}{\color[HTML]{006100} 10000} & \cellcolor[HTML]{C6EFCE}{\color[HTML]{006100} 10000} \\
\hline
\end{tabular}
\caption[Spring: median valid time with clean training data]{Spring: median valid time from cross-validation trials.}\label{spring-vt-table}\end{center}\end{table}

\begin{table}[H]
\begin{center}
\renewcommand{\arraystretch}{0.75}
\begin{tabular}{|c|c|cccccc|}
\hline
 & $r$ \ \textbackslash \ $w_o$ & 0 & $\num{1e-4}$ & $\num{1e-2}$ & $\num{1e-1}$ & 0.5 & 1.0\\ 
\hline
\multirow{5}{*}{$\mathbf{h}_1$} & $\num{1e-12}$ & \cellcolor[HTML]{C6EFCE}{\color[HTML]{006100} 10000} & \cellcolor[HTML]{C6EFCE}{\color[HTML]{006100} 10000} & \cellcolor[HTML]{C6EFCE}{\color[HTML]{006100} 10000} & \cellcolor[HTML]{C6EFCE}{\color[HTML]{006100} 10000} & \cellcolor[HTML]{C6EFCE}{\color[HTML]{006100} 10000} & \cellcolor[HTML]{C6EFCE}{\color[HTML]{006100} 10000} \\
 & $\num{1e-8}$ & \cellcolor[HTML]{C6EFCE}{\color[HTML]{006100} 10000} & \cellcolor[HTML]{C6EFCE}{\color[HTML]{006100} 10000} & \cellcolor[HTML]{C6EFCE}{\color[HTML]{006100} 10000} & 8509 & 6357 & 5925 \\
 & $\num{1e-4}$ & 702 & 702 & 702 & \cellcolor[HTML]{C6EFCE}{\color[HTML]{006100} 703} & 699 & 695 \\
 & $\num{1e-2}$ & 100 & 100 & 100 & 103 & 120 & \cellcolor[HTML]{C6EFCE}{\color[HTML]{006100} 697} \\
 & $\num{1e-1}$ & \cellcolor[HTML]{C6EFCE}{\color[HTML]{006100} 55} & \cellcolor[HTML]{C6EFCE}{\color[HTML]{006100} 55} & \cellcolor[HTML]{C6EFCE}{\color[HTML]{006100} 55} & \cellcolor[HTML]{C6EFCE}{\color[HTML]{006100} 55} & 53 & 52 \\
 \hline
\multirow{5}{*}{$\mathbf{h}_2$} & $\num{1e-12}$ & \cellcolor[HTML]{C6EFCE}{\color[HTML]{006100} 10000} & \cellcolor[HTML]{C6EFCE}{\color[HTML]{006100} 10000} & \cellcolor[HTML]{C6EFCE}{\color[HTML]{006100} 10000} & \cellcolor[HTML]{C6EFCE}{\color[HTML]{006100} 10000} & \cellcolor[HTML]{C6EFCE}{\color[HTML]{006100} 10000} & \cellcolor[HTML]{C6EFCE}{\color[HTML]{006100} 10000} \\
 & $\num{1e-8}$ & \cellcolor[HTML]{C6EFCE}{\color[HTML]{006100} 10000} & \cellcolor[HTML]{C6EFCE}{\color[HTML]{006100} 10000} & \cellcolor[HTML]{C6EFCE}{\color[HTML]{006100} 10000} & \cellcolor[HTML]{C6EFCE}{\color[HTML]{006100} 10000} & \cellcolor[HTML]{C6EFCE}{\color[HTML]{006100} 10000} & \cellcolor[HTML]{C6EFCE}{\color[HTML]{006100} 10000} \\
 & $\num{1e-4}$ & 655 & 655 & 655 & 655 & 657 & \cellcolor[HTML]{C6EFCE}{\color[HTML]{006100} 660} \\
 & $\num{1e-2}$ & 97 & 97 & 98 & 101 & 117 & \cellcolor[HTML]{C6EFCE}{\color[HTML]{006100} 708} \\
 & $\num{1e-1}$ & \cellcolor[HTML]{C6EFCE}{\color[HTML]{006100} 55} & \cellcolor[HTML]{C6EFCE}{\color[HTML]{006100} 55} & \cellcolor[HTML]{C6EFCE}{\color[HTML]{006100} 55} & 54 & 52 & 50 \\
 \hline
\multirow{5}{*}{$\mathbf{h}_3$} & $\num{1e-12}$ & \cellcolor[HTML]{C6EFCE}{\color[HTML]{006100} 10000} & \cellcolor[HTML]{C6EFCE}{\color[HTML]{006100} 10000} & \cellcolor[HTML]{C6EFCE}{\color[HTML]{006100} 10000} & \cellcolor[HTML]{C6EFCE}{\color[HTML]{006100} 10000} & \cellcolor[HTML]{C6EFCE}{\color[HTML]{006100} 10000} & \cellcolor[HTML]{C6EFCE}{\color[HTML]{006100} 10000} \\
 & $\num{1e-8}$ & \cellcolor[HTML]{C6EFCE}{\color[HTML]{006100} 10000} & \cellcolor[HTML]{C6EFCE}{\color[HTML]{006100} 10000} & \cellcolor[HTML]{C6EFCE}{\color[HTML]{006100} 10000} & \cellcolor[HTML]{C6EFCE}{\color[HTML]{006100} 10000} & \cellcolor[HTML]{C6EFCE}{\color[HTML]{006100} 10000} & \cellcolor[HTML]{C6EFCE}{\color[HTML]{006100} 10000} \\
 & $\num{1e-4}$ & 872 & 872 & 872 & 873 & 981 & \cellcolor[HTML]{C6EFCE}{\color[HTML]{006100} 1089} \\
 & $\num{1e-2}$ & 564 & 564 & 568 & \cellcolor[HTML]{C6EFCE}{\color[HTML]{006100} 764} & 597 & 298 \\
 & $\num{1e-1}$ & \cellcolor[HTML]{C6EFCE}{\color[HTML]{006100} 60} & \cellcolor[HTML]{C6EFCE}{\color[HTML]{006100} 60} & \cellcolor[HTML]{C6EFCE}{\color[HTML]{006100} 60} & 59 & 58 & 58 \\
\hline
\end{tabular}
\caption[Lotka-Volterra: median valid time with clean training data]{Lotka-Volterra: median valid time from cross-validation trials.}\label{lv-vt-table}\end{center}\end{table}

\begin{table}[H]
\begin{center}
\renewcommand{\arraystretch}{0.75}
\begin{tabular}{|c|c|cccccc|}
\hline
 & $r$ \ \textbackslash \ $w_o$ & 0 & $\num{1e-4}$ & $\num{1e-2}$ & $\num{1e-1}$ & 0.5 & 1.0\\ 
\hline
\multirow{5}{*}{$\mathbf{h}_1$} & $\num{1e-12}$ & 7890 & 8958 & 8821 & \cellcolor[HTML]{C6EFCE}{\color[HTML]{006100} 10000} & \cellcolor[HTML]{C6EFCE}{\color[HTML]{006100} 10000} & \cellcolor[HTML]{C6EFCE}{\color[HTML]{006100} 10000} \\
 & $\num{1e-8}$ & 5846 & 5849 & 5880 & 7334 & 7411 & \cellcolor[HTML]{C6EFCE}{\color[HTML]{006100} 7498} \\
 & $\num{1e-4}$ & 3904 & 2508 & 3777 & 7342 & \cellcolor[HTML]{C6EFCE}{\color[HTML]{006100} 7846} & 7352 \\
 & $\num{1e-2}$ & 948 & 959 & 2060 & \cellcolor[HTML]{C6EFCE}{\color[HTML]{006100} 4902} & 3754 & 3775 \\
 & $\num{1e-1}$ & 505 & 519 & 996 & 2056 & 3656 & \cellcolor[HTML]{C6EFCE}{\color[HTML]{006100} 4893} \\
\hline
\multirow{5}{*}{$\mathbf{h}_2$} & $\num{1e-12}$ & 7890 & 8958 & 8821 & \cellcolor[HTML]{C6EFCE}{\color[HTML]{006100} 10000} & \cellcolor[HTML]{C6EFCE}{\color[HTML]{006100} 10000} & \cellcolor[HTML]{C6EFCE}{\color[HTML]{006100} 10000} \\
 & $\num{1e-8}$ & 5846 & 5849 & 5880 & 7334 & 7411 & \cellcolor[HTML]{C6EFCE}{\color[HTML]{006100} 7498} \\
 & $\num{1e-4}$ & 3904 & 2508 & 3777 & 7342 & \cellcolor[HTML]{C6EFCE}{\color[HTML]{006100} 7846} & 7352 \\
 & $\num{1e-2}$ & 948 & 959 & 2060 & \cellcolor[HTML]{C6EFCE}{\color[HTML]{006100} 4902} & 3754 & 3775 \\
 & $\num{1e-1}$ & 505 & 519 & 996 & 2056 & 3656 & \cellcolor[HTML]{C6EFCE}{\color[HTML]{006100} 4893} \\
\hline
\multirow{5}{*}{$\mathbf{h}_3$} & $\num{1e-12}$ & 5107 & 7362 & 7955 & 8091 & 8817 & \cellcolor[HTML]{C6EFCE}{\color[HTML]{006100} 8835} \\
 & $\num{1e-8}$ & 2511 & 3711 & 3739 & 3788 & 3783 & \cellcolor[HTML]{C6EFCE}{\color[HTML]{006100} 3803} \\
 & $\num{1e-4}$ & 477 & 513 & \cellcolor[HTML]{C6EFCE}{\color[HTML]{006100} 3833} & 2300 & 2266 & 2293 \\
 & $\num{1e-2}$ & 143 & 147 & 281 & 807 & 2214 & \cellcolor[HTML]{C6EFCE}{\color[HTML]{006100} 3708} \\
 & $\num{1e-1}$ & 87 & 88 & 118 & 216 & 775 & \cellcolor[HTML]{C6EFCE}{\color[HTML]{006100} 815} \\
\hline
\end{tabular}
\caption[Lorenz: median valid time with clean training data]{Lorenz: median valid time from cross-validation trials.}\label{lorenz-vt-table}\end{center}\end{table}

\begin{table}[H]\label{spring-et-table}
\begin{center}
\renewcommand{\arraystretch}{0.77}
\begin{tabular}{|c|c|cccccc|}
\hline
 & $r$ \ \textbackslash \ $w_o$ & 0 & $\num{1e-4}$ & $\num{1e-2}$ & $\num{1e-1}$ & 0.5 & 1.0\\ 
\hline
\multirow{5}{*}{$\mathbf{h}_1$} & $\num{1e-12}$ & $\num{6.5e-10}$ & $\num{6.5e-10}$ & $\num{6.5e-10}$ & $\num{6.5e-10}$ & $\num{6.5e-10}$ & \cellcolor[HTML]{C6EFCE}{\color[HTML]{006100} $\num{6.5e-10}$} \\
 & $\num{1e-8}$ & $\num{6.5e-10}$ & $\num{6.5e-10}$ & $\num{6.5e-10}$ & $\num{6.5e-10}$ & $\num{6.5e-10}$ & \cellcolor[HTML]{C6EFCE}{\color[HTML]{006100} $\num{6.5e-10}$} \\
 & $\num{1e-4}$ & $\num{6.6e-10}$ & $\num{6.6e-10}$ & \cellcolor[HTML]{C6EFCE}{\color[HTML]{006100} $\num{6.5e-10}$} & $\num{6.5e-10}$ & $\num{6.5e-10}$ & $\num{6.5e-10}$ \\
 & $\num{1e-2}$ & $\num{2.8e-9}$ & $\num{2.7e-9}$ & $\num{1.0e-9}$ & \cellcolor[HTML]{C6EFCE}{\color[HTML]{006100} $\num{6.1e-10}$} & $\num{6.4e-10}$ & $\num{6.4e-10}$ \\
 & $\num{1e-1}$ & $\num{2.0e-8}$ & $\num{1.9e-8}$ & $\num{5.2e-9}$ & $\num{9.4e-10}$ & $\num{6.1e-10}$ & \cellcolor[HTML]{C6EFCE}{\color[HTML]{006100} $\num{6.1e-10}$} \\
 \hline
\multirow{5}{*}{$\mathbf{h}_2$} & $\num{1e-12}$ & $\num{1.1e-30}$ & $\num{1.1e-30}$ & $\num{3.7e-31}$ & $\num{2.0e-31}$ & \cellcolor[HTML]{C6EFCE}{\color[HTML]{006100} $\num{2.0e-31}$} & $\num{2.0e-31}$ \\
 & $\num{1e-8}$ & $\num{8.3e-23}$ & $\num{8.1e-23}$ & $\num{1.4e-23}$ & $\num{3.9e-25}$ & $\num{1.8e-26}$ & \cellcolor[HTML]{C6EFCE}{\color[HTML]{006100} $\num{4.6e-27}$} \\
 & $\num{1e-4}$ & $\num{8.3e-15}$ & $\num{8.1e-15}$ & $\num{1.4e-15}$ & $\num{3.9e-17}$ & $\num{1.8e-18}$ & \cellcolor[HTML]{C6EFCE}{\color[HTML]{006100} $\num{4.6e-19}$} \\
 & $\num{1e-2}$ & $\num{8.3e-11}$ & $\num{8.1e-11}$ & $\num{1.4e-11}$ & $\num{3.9e-13}$ & $\num{1.8e-14}$ & \cellcolor[HTML]{C6EFCE}{\color[HTML]{006100} $\num{4.6e-15}$} \\
 & $\num{1e-1}$ & $\num{8.3e-9}$ & $\num{8.0e-9}$ & $\num{1.4e-9}$ & $\num{3.9e-11}$ & $\num{1.8e-12}$ & \cellcolor[HTML]{C6EFCE}{\color[HTML]{006100} $\num{4.6e-13}$} \\
 \hline
\multirow{5}{*}{$\mathbf{h}_3$} & $\num{1e-12}$ & $\num{1.5e-30}$ & $\num{1.4e-30}$ & $\num{4.2e-31}$ & $\num{2.1e-31}$ & \cellcolor[HTML]{C6EFCE}{\color[HTML]{006100} $\num{2.0e-31}$} & $\num{2.0e-31}$ \\
 & $\num{1e-8}$ & $\num{1.2e-22}$ & $\num{1.1e-22}$ & $\num{2.0e-23}$ & $\num{5.3e-25}$ & $\num{2.4e-26}$ & \cellcolor[HTML]{C6EFCE}{\color[HTML]{006100} $\num{6.1e-27}$} \\
 & $\num{1e-4}$ & $\num{1.2e-14}$ & $\num{1.1e-14}$ & $\num{2.0e-15}$ & $\num{5.3e-17}$ & $\num{2.4e-18}$ & \cellcolor[HTML]{C6EFCE}{\color[HTML]{006100} $\num{6.1e-19}$} \\
 & $\num{1e-2}$ & $\num{1.2e-10}$ & $\num{1.1e-10}$ & $\num{2.0e-11}$ & $\num{5.3e-13}$ & $\num{2.4e-14}$ & \cellcolor[HTML]{C6EFCE}{\color[HTML]{006100} $\num{6.1e-15}$} \\
 & $\num{1e-1}$ & $\num{1.2e-8}$ & $\num{1.1e-8}$ & $\num{2.0e-9}$ & $\num{5.3e-11}$ & $\num{2.4e-12}$ & \cellcolor[HTML]{C6EFCE}{\color[HTML]{006100} $\num{6.1e-13}$} \\
\hline
\end{tabular}
\caption[Spring: median discrete energy]{Spring: median discrete energy, $E_h(10,000)$, from cross-validation trials.}\end{center}\end{table}

\begin{table}[H]\label{lv-et-table}
\begin{center}
\renewcommand{\arraystretch}{0.77}
\begin{tabular}{|c|c|cccccc|}
\hline
 & $r$ \ \textbackslash \ $w_o$ & 0 & $\num{1e-4}$ & $\num{1e-2}$ & $\num{1e-1}$ & 0.5 & 1.0\\ 
\hline
\multirow{5}{*}{$\mathbf{h}_1$} & $\num{1e-12}$ & $\num{2.3e-14}$ & $\num{2.3e-14}$ & $\num{1.7e-14}$ & $\num{6.1e-15}$ & $\num{2.1e-15}$ & \cellcolor[HTML]{C6EFCE}{\color[HTML]{006100} $\num{1.2e-15}$} \\
 & $\num{1e-8}$ & $\num{3.0e-10}$ & $\num{3.0e-10}$ & $\num{2.3e-10}$ & $\num{7.9e-11}$ & $\num{2.0e-11}$ & \cellcolor[HTML]{C6EFCE}{\color[HTML]{006100} $\num{1.1e-11}$} \\
 & $\num{1e-4}$ & $\num{1.2e-7}$ & $\num{1.2e-7}$ & $\num{1.1e-7}$ & $\num{7.4e-8}$ & $\num{4.6e-8}$ & \cellcolor[HTML]{C6EFCE}{\color[HTML]{006100} $\num{4.5e-8}$} \\
 & $\num{1e-2}$ & $\num{7.8e-7}$ & $\num{7.8e-7}$ & $\num{7.5e-7}$ & \cellcolor[HTML]{C6EFCE}{\color[HTML]{006100} $\num{6.3e-7}$} & $\num{3.4e-6}$ & $\num{4.6e-5}$ \\
 & $\num{1e-1}$ & $\num{3.9e-1}$ & $\num{3.9e-1}$ & $\num{3.9e-1}$ & $\num{2.8e-1}$ & $\num{2.5e-2}$ & \cellcolor[HTML]{C6EFCE}{\color[HTML]{006100} $\num{1.1e-5}$} \\
 \hline
\multirow{5}{*}{$\mathbf{h}_2$} & $\num{1e-12}$ & $\num{1.2e-15}$ & $\num{1.2e-15}$ & $\num{1.1e-15}$ & $\num{7.7e-16}$ & $\num{4.6e-16}$ & \cellcolor[HTML]{C6EFCE}{\color[HTML]{006100} $\num{3.7e-16}$} \\
 & $\num{1e-8}$ & $\num{6.3e-12}$ & $\num{6.3e-12}$ & $\num{6.0e-12}$ & $\num{4.2e-12}$ & $\num{1.7e-12}$ & \cellcolor[HTML]{C6EFCE}{\color[HTML]{006100} $\num{9.8e-13}$} \\
 & $\num{1e-4}$ & $\num{2.2e-8}$ & $\num{2.2e-8}$ & $\num{2.2e-8}$ & $\num{1.9e-8}$ & $\num{1.1e-8}$ & \cellcolor[HTML]{C6EFCE}{\color[HTML]{006100} $\num{6.6e-9}$} \\
 & $\num{1e-2}$ & $\num{5.3e-7}$ & $\num{5.3e-7}$ & $\num{5.2e-7}$ & $\num{4.6e-7}$ & \cellcolor[HTML]{C6EFCE}{\color[HTML]{006100} $\num{3.8e-7}$} & $\num{1.5e-3}$ \\
 & $\num{1e-1}$ & $\num{3.5e-2}$ & $\num{3.4e-2}$ & $\num{4.2e-2}$ & $\num{6.4e-2}$ & $\num{1.8e-5}$ & \cellcolor[HTML]{C6EFCE}{\color[HTML]{006100} $\num{1.5e-5}$} \\
 \hline
\multirow{5}{*}{$\mathbf{h}_3$} & $\num{1e-12}$ & $\num{5.7e-16}$ & $\num{5.7e-16}$ & $\num{5.5e-16}$ & $\num{4.3e-16}$ & $\num{3.0e-16}$ & \cellcolor[HTML]{C6EFCE}{\color[HTML]{006100} $\num{2.4e-16}$} \\
 & $\num{1e-8}$ & $\num{2.0e-12}$ & $\num{2.0e-12}$ & $\num{1.9e-12}$ & $\num{1.3e-12}$ & $\num{6.0e-13}$ & \cellcolor[HTML]{C6EFCE}{\color[HTML]{006100} $\num{4.1e-13}$} \\
 & $\num{1e-4}$ & $\num{8.5e-9}$ & $\num{8.5e-9}$ & $\num{8.3e-9}$ & $\num{6.8e-9}$ & $\num{3.5e-9}$ & \cellcolor[HTML]{C6EFCE}{\color[HTML]{006100} $\num{2.1e-9}$} \\
 & $\num{1e-2}$ & $\num{3.7e-7}$ & $\num{3.7e-7}$ & $\num{3.6e-7}$ & $\num{3.3e-7}$ & $\num{2.4e-7}$ & \cellcolor[HTML]{C6EFCE}{\color[HTML]{006100} $\num{1.9e-7}$} \\
 & $\num{1e-1}$ & $\num{9.0e-6}$ & $\num{9.0e-6}$ & $\num{9.0e-6}$ & $\num{8.4e-6}$ & $\num{6.3e-6}$ & \cellcolor[HTML]{C6EFCE}{\color[HTML]{006100} $\num{4.4e-6}$}\\
\hline
\end{tabular}
\caption[Lotka-Volterra: median discrete energy]{Lotka-Volterra: median discrete energy, $E_h(10,000)$, from cross-validation trials.}\end{center}\end{table}

\begin{table}[H]\label{lorenz-et-table}
\begin{center}
\renewcommand{\arraystretch}{0.77}
\begin{tabular}{|c|c|cccccc|}
\hline
 & $r$ \ \textbackslash \ $w_o$ & 0 & $\num{1e-4}$ & $\num{1e-2}$ & $\num{1e-1}$ & 0.5 & 1.0\\ 
\hline
\multirow{5}{*}{$\mathbf{h}_1$} & $\num{1e-12}$ & $\num{7.2e-10}$ & $\num{6.9e-10}$ & $\num{5.9e-11}$ & $\num{4.0e-12}$ & $\num{1.7e-12}$ & \cellcolor[HTML]{C6EFCE}{\color[HTML]{006100} $\num{1.3e-12}$} \\
 & $\num{1e-8}$ & $\num{4.5e-10}$ & $\num{4.5e-10}$ & \cellcolor[HTML]{C6EFCE}{\color[HTML]{006100} $\num{3.9e-10}$} & $\num{6.0e-10}$ & $\num{7.8e-10}$ & $\num{8.1e-10}$ \\
 & $\num{1e-4}$ & $\num{1.2e-6}$ & $\num{1.0e-6}$ & $\num{7.4e-8}$ & $\num{3.3e-9}$ & $\num{5.4e-10}$ & \cellcolor[HTML]{C6EFCE}{\color[HTML]{006100} $\num{4.7e-10}$} \\
 & $\num{1e-2}$ & $\num{1.3e-5}$ & $\num{1.4e-5}$ & $\num{1.7e-6}$ & $\num{4.7e-7}$ & $\num{1.4e-7}$ & \cellcolor[HTML]{C6EFCE}{\color[HTML]{006100} $\num{8.1e-8}$} \\
 & $\num{1e-1}$ & $\num{1.2e-4}$ & $\num{9.5e-5}$ & $\num{9.5e-6}$ & $\num{2.0e-6}$ & $\num{6.8e-7}$ & \cellcolor[HTML]{C6EFCE}{\color[HTML]{006100} $\num{4.7e-7}$} \\
 \hline
\multirow{5}{*}{$\mathbf{h}_2$} & $\num{1e-12}$ & $\num{7.2e-10}$ & $\num{6.9e-10}$ & $\num{5.9e-11}$ & $\num{4.0e-12}$ & $\num{1.7e-12}$ & \cellcolor[HTML]{C6EFCE}{\color[HTML]{006100} $\num{1.3e-12}$} \\
 & $\num{1e-8}$ & $\num{4.5e-10}$ & $\num{4.5e-10}$ & \cellcolor[HTML]{C6EFCE}{\color[HTML]{006100} $\num{3.9e-10}$} & $\num{6.0e-10}$ & $\num{7.8e-10}$ & $\num{8.1e-10}$ \\
 & $\num{1e-4}$ & $\num{1.2e-6}$ & $\num{1.0e-6}$ & $\num{7.4e-8}$ & $\num{3.3e-9}$ & $\num{5.4e-10}$ & \cellcolor[HTML]{C6EFCE}{\color[HTML]{006100} $\num{4.7e-10}$} \\
 & $\num{1e-2}$ & $\num{1.3e-5}$ & $\num{1.4e-5}$ & $\num{1.7e-6}$ & $\num{4.7e-7}$ & $\num{1.4e-7}$ & \cellcolor[HTML]{C6EFCE}{\color[HTML]{006100} $\num{8.1e-8}$} \\
 & $\num{1e-1}$ & $\num{1.2e-4}$ & $\num{9.5e-5}$ & $\num{9.5e-6}$ & $\num{2.0e-6}$ & $\num{6.8e-7}$ & \cellcolor[HTML]{C6EFCE}{\color[HTML]{006100} $\num{4.7e-7}$} \\
 \hline
\multirow{5}{*}{$\mathbf{h}_3$} & $\num{1e-12}$ & $\num{9.7e-9}$ & $\num{7.1e-9}$ & $\num{6.2e-10}$ & $\num{2.9e-10}$ & \cellcolor[HTML]{C6EFCE}{\color[HTML]{006100} $\num{1.5e-10}$} & $\num{1.5e-10}$ \\
 & $\num{1e-8}$ & $\num{2.0e-6}$ & $\num{2.0e-6}$ & $\num{1.6e-6}$ & $\num{7.8e-7}$ & $\num{2.4e-7}$ & \cellcolor[HTML]{C6EFCE}{\color[HTML]{006100} $\num{1.2e-7}$} \\
 & $\num{1e-4}$ & $\num{2.7e-4}$ & $\num{1.8e-4}$ & $\num{7.8e-6}$ & \cellcolor[HTML]{C6EFCE}{\color[HTML]{006100} $\num{3.9e-6}$} & $\num{4.3e-6}$ & $\num{4.3e-6}$ \\
 & $\num{1e-2}$ & $\num{2.4e-3}$ & $\num{5.2e-3}$ & $\num{1.7e-3}$ & $\num{7.7e-5}$ & $\num{1.1e-5}$ & \cellcolor[HTML]{C6EFCE}{\color[HTML]{006100} $\num{6.8e-6}$} \\
 & $\num{1e-1}$ & $\num{2.5e-2}$ & $\num{1.7e-2}$ & $\num{4.8e-3}$ & $\num{1.7e-3}$ & $\num{3.5e-4}$ & \cellcolor[HTML]{C6EFCE}{\color[HTML]{006100} $\num{9.3e-5}$}\\
\hline
\end{tabular}
\caption[Lorenz: median discrete energy]{Lorenz: median discrete energy, $E_h(10,000)$, from cross-validation trials.}\end{center}\end{table}

Tables \ref{lv-et-table} and \ref{lorenz-et-table} generally corroborate this result for the other test problems, although the magnitude of median discrete energy improvement varies between problems and state functions. In the case of the Lotka-Volterra problem, we are pleased to see orderly and favorable results for physics-informed training under the discrete energy metric, despite apparently poor performance under the valid time metric. While median discrete energy does increase with regularization, meaningful physics-informed training is able to salvage respectable ODE adherence error not worse than the order of $\num{1e-5}$ for the superior models. These findings suggest that while the recursive prediction breaches the $M=\num{1e-4}$ validity threshold early, the prediction displays good adherence to the governing ODE. In the case of the chaotic Lorenz problem, Table \ref{lorenz-et-table} provides firm confirmation that physics-informed training yields predictions which adhere to the ODE well, strengthening the data-driven findings from Table \ref{lorenz-vt-table}. Since the reference data for this problem is inherently untrustworthy, the physics-informed training adds trustworthy structure to the model. Altogether, these results conform with our findings using only a data-driven metric.

\section{Conclusions}\label{sec-conclusions}

In this article, we began by reviewing multiple classes of numerical methods for ordinary differential equations spanning familiar numerical integration schemes, more complex recurrent neural networks, and novel reservoir computers. Recent work and theory was considered regarding neural network solutions to partial differential equations. The framework was applied to NVAR to derive the companion model with shared parameters that enforces the structure of the right-hand side of the underlying differential equation. We identify this secondary model as the physics-informed NVAR.

The effectiveness of piNVAR was tested in multiple numerical experiments. We developed a central cross-validation testing routine to examine model predictive performance and generalization for multiple ODEs, which ranged from linear to nonlinear and chaotic, using two metrics, data-driven valid time and physics-informed discrete energy. For multiple test model structures with no incorporation of the right-hand side of the differential equation into the state function, we showed that models trained with large ODE-based training weights generally exhibit superior predictive performance for most choices of regularization. Lastly, we found ODE training also substantially improved model performance under the physics-informed discrete energy metric, and confirmed the validity of the model predictions through this reference-agnostic lens.

Future work includes the analysis of higher-order polynomial state functions, varied parameters for the support coefficient and input transformation components of the state functions, and the comparison of multiple model sizes.

\section*{Code availability}

Source code for all results shown in this article can be found at \url{https://github.com/samuelhocking/pinvar-sisc-2024}. All code is released under an open source MIT License.

\section*{Acknowledgments}

This research was supported, in part, by a seed grant from Tufts University (Data Intensive Studies Center) and the National Science Foundation (NSF 2008276).

\newpage
\bibliographystyle{siamplain}

\end{document}